\documentclass[10pt, a4paper,leqno]{amsart}

\usepackage{amsmath,amsthm,amscd,amssymb,amsfonts, amsbsy, mathtools}
\usepackage{latexsym}
\usepackage{exscale}
\usepackage{mathrsfs}
\usepackage{euscript}
\usepackage{enumerate}

\usepackage{pgf}
\usepackage[utf8]{inputenc}
\usepackage[OT2,T1]{fontenc}
\usepackage{lmodern}

\usepackage[colorlinks=true, pdfstartview=FitV, linkcolor=blue, citecolor=red, urlcolor=blue]{hyperref}
\usepackage[doipre={DOI:~}]{uri}

\theoremstyle{plain}

\newtheorem*{Cuculescutheo}{Cuculescu's construction \cite{cuculescu1971}}

\newtheorem*{theo*}{Theorem}

\numberwithin{equation}{section}

\newtheorem{innercustomthm}{Theorem}
\newenvironment{customthm}[1]
  {\renewcommand\theinnercustomthm{#1}\innercustomthm}
  {\endinnercustomthm}

\theoremstyle{remark}
\newtheorem{remark}[equation]{Remark}

\theoremstyle{definition}

\newtheorem*{def*}{Definition}

\newcommand{\N}{\mathbb{N}}
\newcommand{\Z}{\mathbb{Z}}
\newcommand{\R}{\mathbb{R}}

\DeclareSymbolFont{cyrletters}{OT2}{wncyr}{m}{n}
\DeclareMathSymbol{\Sha}{\mathalpha}{cyrletters}{"58}

\def\esssup{\mathop{\mathrm{ess \ sup}}}

\newcommand{\minus}{\scalebox{0.75}[1.0]{$-$}}

\begin{document}

\title[Operator-valued dyadic harmonic analysis]{Operator-valued dyadic harmonic analysis beyond doubling measures}

\author[J.M. Conde-Alonso]{Jos\'e M. Conde-Alonso}
\address{Jos\'e M. Conde-Alonso\\
Instituto de Ciencias Matem\'{a}ticas CSIC-UAM-UC3M-UCM\\
Consejo Superior de Investigaciones Cient\'{\i}ficas\\
C/ Nicol\'{a}s Ca\-brera, 13-15\\
E-28049 Madrid, Spain}

\email{jose.conde@icmat.es}

\

\author[L.D. L\'{o}pez-S\'{a}nchez]{Luis Daniel L\'{o}pez-S\'{a}nchez}
\address{Luis Daniel L\'{o}pez-S\'{a}nchez\\
Instituto de Ciencias Matem\'{a}ticas CSIC-UAM-UC3M-UCM\\
Consejo Superior de Investigaciones Cient\'{\i}ficas\\
C/ Nicol\'{a}s Ca\-brera, 13-15\\
E-28049 Madrid, Spain}

\email{luisd.lopez@icmat.es}

\

\date{\today}

\thanks{Partially supported by the European Research Council ERC StG-256997-CZOSQP, the Spanish grant MTM2010-16518 and by ICMAT Severo Ochoa Grant SEV-2011-0087 (Spain)}




\begin{abstract}
We obtain a complete characterization of the weak-type $(1,1)$ for Haar shift operators in terms of generalized Haar systems adapted to a Borel measure $\mu$ in the operator-valued setting. The main technical tool in our method is a noncommutative Calder\'on-Zygmund decomposition valid for arbitrary Borel measures.
\end{abstract}
\maketitle

\section{Introduction}
\label{section:Intro}

We say that $\Phi=\{\phi_Q\}_{Q\in\mathscr{D}}$ is a \emph{generalized Haar system} in $\R^d$ adapted to a locally finite Borel measure $\mu$ and a dyadic lattice $\mathscr{D}$ if the following conditions hold:
\begin{enumerate}[\normalfont(a)]\leftmargin=.8cm
\labelwidth=.8cm\itemsep=0.2cm\topsep=.1cm

\item For every $Q\in\mathscr{D}$, $\mathrm{supp}(\phi_Q)\subset Q$.

\item If $Q'$, $Q\in\mathscr{D}$ and $Q'\subsetneq Q$, then $\phi_Q$ is constant on $Q'$.

\item \label{item:Haarcan} For every $Q\in\mathscr{D}$, $\displaystyle\int_{\R^d} \phi_Q \,d\mu = 0.$

\item \label{Haar:propd} For every $Q\in\mathscr{D}$, either $\|\phi_Q\|_{L^2(\mu)}=1$ or $\phi_Q\equiv 0$ and $\mu(Q)=0$.
\end{enumerate}
If the vanishing integral condition \eqref{item:Haarcan} is not imposed, the Haar system is said to be \emph{non-cancellative}. Let $\Phi = \{\phi_Q\}_{Q \in \mathscr{D}}$ and $\Psi = \{\psi_Q\}_{Q \in \mathscr{D}}$ be two non-necessarily cancellative generalized Haar systems in $\R^d$. A Haar shift operator of complexity $(r,s) \in \N \times \N$ is an operator of the form
\begin{equation}\label{eq:cSha}
\Sha_{r,s} f (x) = \sum_{Q \in \mathscr{D}} \sum_{\begin{subarray}{c} R \in \mathscr{D}_r(Q) \\ S \in \mathscr{D}_s(Q) \end{subarray}} \alpha_{R,S}^Q \langle f, \phi_R \rangle \psi_S (x), \quad \textrm{with} \quad \sup_{Q,R,S} |\alpha^Q_{R,S}| < \infty;
\end{equation}
where $\langle f,g\rangle = \int_{\R^d} fg \, d\mu$ and $\mathscr{D}_k(Q)$, $k \in \N$, denotes the family of $k$-dyadic descendants of $Q$: the partition of $Q$ into subcubes $R \in \mathscr{D}$ of side-length $\ell(R) = 2^{-k}\ell(Q)$. Several objects in dyadic harmonic analysis have the general form \eqref{eq:cSha}, including Haar multipliers, dyadic paraproducts, the dyadic model of the Hilbert transform and their adjoints. Haar shift operators have served as important tools in the study of many different problems in harmonic analysis since the form \eqref{eq:cSha} is a fruitful source of  models of Calder\'on-Zygmund operators. In particular, in the case where $\mu$ is the Lebesgue measure, Calder\'on-Zygmund operators can be expressed as weak limits of certain averages of cancellative Haar shift operators and paraproducts \cite{hytonen2012} and are pointwise dominated by positive dyadic operators, which are Haar shift o\-pe\-ra\-tors relative to non-cancellative Haar systems \cite{condealoso-rey2014}.

The boundedness behavior of Haar shift operators with respect to arbitrary locally finite Borel measures in the commutative setting was studied in \cite{lopezsanchez-martell-mparcet2014}. There the authors characterize the weak-type $(1,1)$ of such operators. In this note we extend the scope of this result to the setting of semicommutative $L_p$ spaces. The main technique that we will use in our approach is a generalization of the Calder\'on-Zygmund decomposition introduced in \cite{lopezsanchez-martell-mparcet2014} which is valid for operator-valued functions, in the spirit of the Calder\'{o}n-Zygmund decomposition constructed in \cite{parcet2009}.

We will work in the following framework: consider a pair $(\mathcal{M}, \nu)$ where $\mathcal{M}$ is a von Neumann algebra and $\nu$ is a normal semifinite faithful trace on $\mathcal{M}$ and let $\mu$ be a locally finite Borel measure on $\R^d$. Let $\mathcal{A}_B$  be the algebra of essentially bounded $\mathcal{M}$-valued functions
\begin{equation*}
\label{algebra}
\mathcal{A}_B = \Biggl\{ f: \R^d \to \mathcal{M} \, \big| \, f \ \mbox{strongly measurable s.t.} \ \esssup_{x \in \R^d} \|f(x)\|_{\mathcal{M}} < \infty \Biggr\}
\end{equation*}
equipped with the \emph{n.s.f.} trace $\tau(f) = \int_{\R^d} \nu(f) \,d\mu$. The weak-operator closure $\mathcal{A}$ of $\mathcal{A}_B$ is a von Neumann algebra isomorphic to $L_{\infty}(\R^d,\mu) \overline{\otimes} \mathcal{M}$. Given a rearrangement invariant quasi-Banach function space $X$, let us write $X(\mathcal{M})$ and $X(\mathcal{A})$ for their associated noncommutative symmetric spaces. In particular $L_p(\mathcal{M})$ and $L_p(\mathcal{A})$ denote the noncommutative $L_p$ spaces associated to the pairs $(\mathcal{M},\nu)$ and $(\mathcal{A},\tau)$. It can be readily seen that for $1 \leq p < \infty$ the noncommutative $L_p$ space $L_p(\mathcal{A})$ is isometric to the Bochner $L_p$ space $L_p(\R^d, \mu ; L_p(\mathcal{M}))$. The lattices of projections are denoted by $\mathcal{P}(\mathcal{M})$ and $\mathcal{P}(\mathcal{A})$, while $1_\mathcal{M}$ and $1_\mathcal{A}$ stand for the unit elements and $\mathcal{M}'$ and $\mathcal{A}'$ stand for their respective commutants. For a more detailed discussion on noncommutative $L_p$ spaces we refer to \cite{mei-parcet2009} and references therein. The reader unfamiliar with the theory of noncommutative $L_p$ spaces may think of $\mathcal{M}$ as the algebra $\mathcal{B}(\ell_2^n)$ of $n \times n$ matrices equipped with the standard trace $Tr$, thereby recovering the classical Schatten $p$-classes. The reader should take into account that, with this setting in mind, we provide estimates uniform on $n$.

Before stating our results let us introduce some notation first. By $(\mathsf{E}_k)_{k \in \Z}$ we will denote the family of conditional expectations associated to $\mathscr{D}_k$ --- the dyadic cubes $Q$ of side-length $\ell(Q) = 2^{-k}$ --- and write $\mathsf{D}_k$ for the corresponding martingale difference operators.  The tensor product $\mathsf{E}_k \otimes id_\mathcal{M}$ acting on $\mathcal{A}$ will also be denoted by  $\mathsf{E}_k$, which yields a filtration $(\mathcal{A}_k)_{k \in \Z}$ on $\mathcal{A}$. We thus have that
\begin{align*}
\mathsf{E}_k (f) & =: f_k =  \sum_{Q \in \mathscr{D}_k} \langle f \rangle_Q  1_Q, 
\\ 
\mathsf{D}_k (f) & =: df_k = \sum_{Q \in \mathscr{D}_k} \big( \langle f \rangle_Q - \langle f \rangle_{\widehat{Q}} \big) 1_Q,
\end{align*} 
which correspond to projections to the class of operators constant at scale $\mathscr{D}_k$. Here $1_Q$ denotes the characteristic function of $Q$, $\langle f \rangle_Q = \mu(Q)^{-1} \int_Q f \,d\mu$ and $\widehat{Q}$ is the dyadic parent of $Q$: the only dyadic cube that contains $Q$ with twice its side-length.

We will construct the Calder\'on-Zygmund decomposition for functions in the class
\begin{equation*}
\label{denseclass}
\mathcal{A}_{+,K} = \{f: \R^d \to \mathcal{M} \, | \, f \geq 0, \; \mathrm{supp}_{\R^d} (f) \; \textrm{is compact} \},
\end{equation*}
whose span is dense in $L_1(\mathcal{A})$. Here $\mathrm{supp}_{\R^d} (f)$ stands for the support of $f$ as an operator-valued function, as opposed to its support projection as an element of a von Neumann algebra. As the Calder\'on-Zygmund decomposition introduced in \cite{parcet2009} --- which is suitable for the Lebesgue measure and doubling measures --- the Calder\'on-Zygmund decomposition here presented is comprised of diagonal and off-diagonal terms, reflecting the lack of commutativity in the operator-valued framework. Taking $i \vee j = \max\{i,j\}$ and $i \wedge j = \min \{i,j\}$ for $i, j \in \Z$ we have:

\begin{customthm}{A}\label{thm:CZd} Let $f \in \mathcal{A}_{+,K}$ and let $\lambda > 0$. Then there exist a family of pairwise disjoint projections $(p_k)_{k\in \Z}$ adapted to $(\mathcal{A}_k)_{k \in \Z}$ and a projection $q:= 1_{\mathcal{A}} - \sum_k p_k  \in \mathcal{P}(\mathcal{A})$ such that $f$ can be decomposed as $f = g + b + \beta$, where each term has a diagonal and an off-diagonal part given by
\begin{list}{$\bullet$}{\leftmargin=0.4cm\labelwidth=.4cm  \itemsep=0.2cm}
	\item $g =  g_{\Delta} + g_{\mathrm{off}}$, where
		\begin{align*}
		g_{\Delta} &= q f q + \sum_{k \in \Z} \mathsf{E}_{k-1}\left(p_{k} f_k p_{k}\right),  
		\\ 
		g_{\mathrm{off}} &= (1_\mathcal{A} -q)fq+ q f (1_\mathcal{A} -q) + \sum_{i\not=j} \mathsf{E}_{i\vee j-1}\left(p_{i} f_{i\vee j} p_{j}\right);
		\end{align*}
	
	\item $b = b_{\Delta} + b_{\mathrm{off}}$, where
		\begin{equation*}
		b_{\Delta} = \sum_{k \in \Z} p_{k}(f-f_{k})p_{k}\,,
		\qquad
		b_{\mathrm{off}} = \sum_{i\not=j} p_{i}(f-f_{i\vee j})p_{j}; 
		\end{equation*}	

	\item $\beta = \beta_{\Delta} + \beta_{\mathrm{off}}$, where
		\begin{equation*}
		\beta_{\Delta} = \sum_{k \in \Z} \mathsf{D}_{k} (p_{k} f_k p_{k}),
		\qquad
		 \beta_{\mathrm{off}} = \sum_{i\not=j} \mathsf{D}_{i\vee j}\left(p_{i} f_{i \vee j} p_{j}\right).
		\end{equation*}
\end{list}
The diagonal terms satisfy the classical properties
\begin{enumerate}[\normalfont(a)]\leftmargin=.8cm \labelwidth=.8cm\itemsep=0.2cm\topsep=.1cm
	\item \label{item:LdosCZd} $g_{\Delta} \in L_1(\mathcal{A}) \cap L_2(\mathcal{A})$ with
		\begin{equation*}
		\|g_{\Delta}\|_{L_1(\mathcal{A})} = \|f\|_{L_1(\mathcal{A})}, \quad \|g_{\Delta}\|_{L_2(\mathcal{A})}^2 \leq 39 \lambda\|f\|_{L_1(\mathcal{A})};
		\end{equation*}
	\item \label{item:badbCZd} $b_{\Delta}  = \sum_{k \in \Z} b_k$, with $\int_{\R^d} b_k \,d\mu = 0$ and satisfies the estimate
		\begin{equation*}
		\|b_{\Delta}\|_{L_1(\mathcal{A})} =  \sum_{k \in \Z} \|b_k\|_{L_1(\mathcal{A})} \leq 2\|f\|_{L_1(\mathcal{A})};
		\end{equation*}	
	\item \label{item:badbetaCZd} $\beta_{\Delta}  = \sum_{k \in \Z} \beta_k$, with each $\beta_k$ a $k$-th martingale difference, and is such that
		\begin{equation*}
		\|\beta_{\Delta}\|_{L_1(\mathcal{A})} \leq \sum_{k \in \Z}  \|\beta_k\|_{L_1(\mathcal{A})}  \leq 2 \|f\|_{L_1(\mathcal{A})}.
		\end{equation*}
\hspace{-8 mm} 
The off-diagonal terms are such that
	\item \label{item:LdosCZoff} $g_{\mathrm{off}}$ decomposes as $g_{\mathrm{off}}= \sum_{k \in \Z, h\geq 1} g_{k,h}$, where $g_{k,h}$ is the $(k + h)$-th martingale difference $g_{k,h} = \mathsf{D}_{k + h}(p_kf_{k+h}q_{k+h} + q_{k+h}f_{k+h}p_k)$, and satisfies the estimate
			\begin{equation*}
			\sup_{h \geq 1} \sum_{k \in \Z} \|g_{k,h}\|_{L_2(\mathcal{A})}^2 \leq 16\lambda\|f\|_{L_1(\mathcal{A})};
			\end{equation*}
	\item \label{item:badbCZoff} $b_{\mathrm{off}}  = \sum_{k \in \Z,h \geq 1} b_{k,h}$, where $b_{k,h} = p_k(f - f_{k+h})p_{k+h} + p_{k+h}(f - f_{k+h})p_k$, $\int_{\R^d} b_{k,h} \, d\mu = 0$ and
		\begin{equation*}
		\sum_{k \in \Z} \|b_{k,h}\|_{L_1(\mathcal{A})} \leq 8 (h + 1) \|f\|_{L_1(\mathcal{A})};
		\end{equation*}	
	\item \label{item:badbetaCZoff} $\beta_{\mathrm{off}}  = \sum_{k \in \Z,h \geq 1} \beta_{k,h}$, where $\beta_{k,h} = \mathsf{D}_{k+h}(p_{k}f_{k+h}p_{k+h} + p_{k+h}f_{k+h}p_{k})$ and
		\begin{equation*}
		\sum_{k \in \Z} \|\beta_{k,h}\|_{L_1(\mathcal{A})} \leq 8(h + 1) \|f\|_{L_1(\mathcal{A})}.
		\end{equation*}
\end{enumerate}
\end{customthm}

 Observe that the diagonal terms satisfy estimates similar to those of their commutative counterparts found in \cite{lopezsanchez-martell-mparcet2014}. However, in contrast to the classical setting, there are additional difficulties in proving the estimates even for diagonal terms due to the noncommutativity of $\mathcal{A}$. In particular, the estimates of $g_{\Delta}$ are proved in a different way and only hold for $p \leq 2$. In addition, the fact that $\mu$ is allowed to be nondoubling brings other difficulties not present in \cite{parcet2009}. On the other hand, at first glance the off-diagonal estimates in \eqref{item:LdosCZoff}, \eqref{item:badbCZoff} and \eqref{item:badbetaCZoff} seem to be insufficient, since they are weaker than the expected ones: $\|g_{\mathrm{off}}\|_{L_2(\mathcal{A})} \lesssim \lambda \|f\|_{L_1(\mathcal{A})}$, $\sum_{k,h}\|b_{k,h}\|_{L_1(\mathcal{A})} \lesssim \|f\|_{L_1(\mathcal{A})}$ and $\sum_{k,h}\|\beta_{k,h}\|_{L_1(\mathcal{A})} \lesssim \|f\|_{L_1(\mathcal{A})}$. Moreover, estimates of this nature seem to fail as hinted in \cite{parcet2009}. However, the estimates at hand will prove to be sufficient for our purposes as the operators under consideration are localized in a sense stronger than in \cite{mei-parcet2009, parcet2009}. In that respect, one can think of our result as a partial answer to the question posed in \cite{mei-parcet2009} about the existence of a Littlewood-Paley theory for nondoubling measures in the semicommutative context.


Let $\Phi = \{\phi_Q\}_{Q \in \mathscr{D}}$ and $\Psi = \{\psi_Q\}_{Q \in \mathscr{D}}$ be two non-necessarily cancellative generalized Haar systems. A \emph{commuting Haar shift operator} is an $L_2(\mathcal{A})$ bounded operator of the form
\begin{equation}\label{eq:ncHaar}
\Sha_{r,s}f(x)  = \sum_{Q \in \mathscr{D} } \sum_{\begin{subarray}{c} R \in \mathscr{D}_r(Q) \\ S \in \mathscr{D}_s(Q) \end{subarray}} \alpha_{R,S}^Q \langle f, \phi_R \rangle \psi_S(x), \quad \sup_{Q,R,S} \|\alpha_{R,S}^Q\|_{\mathcal{M}} < \infty,
\end{equation}
where the symbols $\alpha_{R,S}^Q$ lie in $\mathcal{M} \cap \mathcal{M}'$, the center of $\mathcal{M}$. Notice that in this definition the pairing $\langle f, g \rangle = \int_{\R^d} fg \,d\mu$ is in fact a partial trace and whence operator-valued. Our second result determines conditions for which the weak-type $(1,1)$ for these operators hold.

\begin{customthm}{B}
\label{thm:HS}
Let $\Sha_{r,s}$ be given as in \eqref{eq:ncHaar}. Assume that $\Sha_{r,s}$ satisfies the restricted local vector-valued $L_2$ estimate
\begin{equation}
\label{eq:localv-v}
\int_{\R^d} \|\Sha_{r,s}^{Q_0}(1_{Q_0})(x)\|^2_{\mathcal{M}} \, d\mu(x) \leq C \mu(Q_0),
\end{equation}
uniformly over $Q_0 \in \mathscr{D}$. Here
\begin{equation*}
\Sha_{r,s}^{Q_0}f(x)  := \sum_{Q \in \mathscr{D}(Q_0) } \sum_{\begin{subarray}{c} R \in \mathscr{D}_r(Q) \\ S \in \mathscr{D}_s(Q) \end{subarray}} \alpha_{R,S}^Q  \langle f,\phi_R\rangle \psi_S(x),
\end{equation*}
where $\mathscr{D}(Q)$ denotes the family of all dyadic subcubes of $Q$ including $Q$ itself. If
\begin{equation}\label{eq:m-funct}
\Xi(\Phi, \Psi, r,s) := \sup_{Q \in \mathscr{D}} \{ \|\phi_R\|_{L_{\infty}(\mu)} \|\psi_S\|_{L_1(\mu)} \,:\, R \in \mathscr{D}_r(Q), S \in \mathscr{D}_s(Q)\} < \infty.
\end{equation}
then $\Sha_{r,s}$ maps $L_1(\mathcal{A})$ continuously into $L_{1,\infty}(\mathcal{A})$.
\end{customthm}

\begin{remark}
A testing argument with simple functions is used in \cite{lopezsanchez-martell-mparcet2014} to show that the condition \eqref{eq:m-funct} is also necessary when the symbols are all nonzero. One can show that this is also the case in the present setting by following similar ideas, and hence they will not be repeated here.	
\end{remark}

\begin{remark}
As in the commutative case, if the Haar systems $\Phi = \{\phi_Q\}_{Q \in \mathscr{D}}$ and $\Psi = \{\psi_Q\}_{Q \in \mathscr{D}}$ are cancellative, orthogonality arguments may be used to verify that the condition \eqref{eq:localv-v} and the $L^2$ boundedness of $\Sha_{r,s}$ are satisfied.
\end{remark}

The condition \eqref{eq:m-funct} may be interpreted as certain restriction on the measure $\mu$ in terms of its degeneracy over generations of dyadic cubes. The resulting class of measures depends strongly on the Haar shift operator in question. For some operators the associated class of measures is shown to be strictly bigger than the doubling class, but nevertheless disjoint from the class of measures of polynomial growth, for which non-standard Calder\'{o}n-Zygmund theories are available; see \cite{lopezsanchez-martell-mparcet2014} and references therein.


\section{The Calder\'on-Zygmund decomposition}
\label{section:CZdec}

This section is devoted to the proof of Theorem \ref{thm:CZd}. First, some reductions are in order. For simplicity we will assume that $\mu(\R^d)=\infty$ and that the dyadic lattice $\mathscr{D}$ has no quadrants. Namely, that $\mathscr{D}$ is such that for every compact $K$ there exists $Q \in \mathscr{D}$ with $K \subset Q$. These assumptions can be removed arguing as in \cite{lopezsanchez-martell-mparcet2014}. However, we find the second assumption very natural since --- in a probabilistic sense --- almost all dyadic lattices satisfy it. Also, as argued  in \cite{lopezsanchez-martell-mparcet2014}, we are confident that our results also hold in the context of geometrically doubling metric spaces. From the previous assumptions, it can be seen that for a fixed $f\in \mathcal{A}_{+,K}$ and $\lambda > 0$ there exists $m_\lambda(f) \in \Z$ such that $f_k \leq \lambda 1_{\mathcal{A}}$ for all $k \leq m_\lambda(f)$ (see \cite{parcet2009}). Without loss of generality, we may also assume that $f$ has only finite non-vanishing martingale differences.

\begin{remark}
To ease notation, we will use the normalization $m_\lambda(f) = 0$. It is safe to assume so since in the proofs of Theorems \ref{thm:CZd} and \ref{thm:HS} both $f \in \mathcal{A}_{+,K}$ and $\lambda > 0$ will remain fixed, but otherwise arbitrary.
\end{remark}

We start with the construction of the projections $(p_k)_{k \in \Z}$ and $q$ of Theorem \ref{thm:CZd}. To that end we will use the so-called Cuculescu's construction. Here we state it in the precise form that we will use, although the construction can be done in any semifinite von Neumann algebra.

\begin{Cuculescutheo}
\label{thm:Cucu}
Let $f \in \mathcal{A}_{+,K}$ and consider the associated positive martingale $(f_k)_{k \in \Z}$ relative to the dyadic filtration $(\mathcal{A}_k)_{k \in \Z}$. Given $\lambda > 0$, the decreasing sequence of projections $(q_{k})_{k \in \Z}$ defined recursively by $q_k = 1_{\mathcal{A}}$ for $k \leq 0$ and
\begin{equation*}
q_{k} = q_k(f,\lambda) :=  1_{(0,\lambda]} \bigl(q_{k-1} f_k q_{k-1} \bigr)
\end{equation*}
is such that
\begin{enumerate}[\normalfont(a)]\leftmargin=.8cm
\labelwidth=.8cm\itemsep=0.2cm\topsep=.1cm
\item $q_{k}$ is a projection in $\mathcal{A}_k$.
\item \label{item:Cucomm} $q_{k}$ commutes with $q_{k-1} f_k q_{k-1}$.
\item \label{item:Cuest} $q_{k} f_k q_{k} \leq \lambda q_{k}$.
\item \label{item:Cucuest} $q = \bigwedge_k q_{k}$ satisfies
\begin{equation*}
\|q f_k q \|_\mathcal{A} \leq \lambda \,\, \textrm{for all} \,\, k \geq 1\quad \textrm{and} \quad \tau(1_{\mathcal{A}} - q) \leq \frac{1}{\lambda}\|f\|_{L_1(\mathcal{A})}.
\end{equation*}
\end{enumerate}
Define the sequence $(p_k)_{k \geq 1}$ of pairwise disjoint projections by
\begin{equation*}
p_{k} = q_{k-1} - q_{k}.
\end{equation*}
In particular
\begin{equation*}
\sum_{k \geq 1} p_k = 1_\mathcal{A} - q
\end{equation*}
and also $p_{k}f_kp_{k} \geq \lambda p_{k}$.
\end{Cuculescutheo}

\begin{remark}
\label{remark:locCucu}
Since the projection $q_k$ is $\mathscr{D}_k$-measurable, we have the following useful expression
\begin{equation*}
q_k = \sum_{Q \in \mathscr{D}_k} q_Q \otimes 1_Q,
\end{equation*}
where $q_Q = q_Q(f,Q)$ are projections in $\mathcal{M}$ defined by
\begin{equation*}
q_Q = \begin{dcases}
1_\mathcal{M} &\mbox{if } k < 0 \\
1_{(0,\lambda]} ( q_{\widehat{Q}} \, \langle f \rangle_Q \, q_{\widehat{Q}} ) &\mbox{if } k \geq 0.
\end{dcases}
\end{equation*}
As in Cuculescu's construction, these projections satisfy
\begin{enumerate}[\normalfont(a)]\leftmargin=.8cm
\labelwidth=.8cm\itemsep=0.2cm\topsep=.1cm
\item\label{item:locCucusub} $q_Q \leq q_{\widehat{Q}}$.
\item $q_Q$ commutes with $q_{\widehat{Q}} \, \langle f \rangle_Q \, q_{\widehat{Q}}$.
\item $q_Q \, \langle f \rangle_Q \, q_Q \leq \lambda q_Q$.
\end{enumerate}
One then can express the projections $p_k$ as
\begin{equation}
\label{eq:CZcubes}
p_{k} = \sum_{Q \in \mathscr{D}_k} (q_{\widehat{Q}} - q_Q)1_Q =: \sum_{Q \in \mathscr{D}_k} p_Q \otimes 1_Q,
\end{equation}
and we analogously have that $p_Q \in \mathcal{P}(\mathcal{M})$ is such that $p_Q \langle f \rangle_Q \,p_Q \geq \lambda p_Q$. As detailed in \cite{parcet2009}, one could interpret the projections $p_k$ as the union dyadic cubes of side-length $2^{-k}$ into which the classical level set $\Omega_\lambda = \{\sup_{k}f_k > \lambda\}$ is decomposed as in \cite{lopezsanchez-martell-mparcet2014}. One can thus view $q$ as the complementary set of $\Omega_\lambda$.
\end{remark}


\subsection*{Proof of Theorem \ref{thm:CZd}}

By construction $f = g + b + \beta$. We now turn to the estimates of the diagonal part. For the $L_1$ estimate of $g_{\Delta}$ observe that by the tracial property
\begin{multline*}
\|g_{\Delta}\|_{L_1(\mathcal{A})} = \tau(fq) + \sum_{k \geq 1} \tau(\mathsf{E}_{k-1}(p_{k}f_kp_{k})) = \tau(fq) + \tau(f(1_\mathcal{A} - q))=  \|f\|_{L_1(\mathcal{A})},
\end{multline*}
since $\mathsf{E}_{k}$ preserves the trace. The proof of the $L_2$ estimate of $g_{\Delta}$ is a bit more involved since $\mu$ is not necessarily doubling. Also, the lack of commutativity of $\mathcal{M}$ prevents us from following the argument that appeared in \cite{lopezsanchez-martell-mparcet2014}. However, standard arguments in noncommutative martingale theory apply. First notice that since $q_k$ commutes with $q_{k-1}f_{k}q_{k-1}$, 
\begin{align*}
\mathsf{E}_{k-1}(p_{k}f_kp_{k}) = q_{k-1}f_{k-1}q_{k-1} - \mathsf{E}_{k-1}(q_{k}f_kq_{k}).
\end{align*}
Thus,
\begin{align*}
\Bigg\|\sum_{k \geq 1} \mathsf{E}_{k-1}(p_{k}f_kp_{k}) \Bigg\|_{L_2(\mathcal{A})}^2 &\leq 2 \Biggl(\Bigg\|\sum_{k \geq 1} q_{k}f_kq_{k} - \mathsf{E}_{k-1}(q_{k}f_kq_{k}) \Bigg\|_{L_2(\mathcal{A})}^2 
\\
&\qquad\quad + \Bigg\|\sum_{k \geq 1} q_{k}f_kq_{k} - q_{k-1}f_{k-1}q_{k-1} \Bigg\|_{L_2(\mathcal{A})}^2 \Biggr) 
\\
&= 2(I + II).
\end{align*}
As it is proved in \cite[Lemma 3.4]{randrianantoanina2002}, we have that
\begin{multline*}
\|q_{k}f_kq_{k} - \mathsf{E}_{k-1}(q_{k}f_kq_{k})\|_{L_2(\mathcal{A})}^2 \leq 2\big(\|q_{k}f_kq_{k}\|_{L_2(\mathcal{A})}^2 -  \|q_{k-1}f_{k-1}q_{k-1}\|_{L_2(\mathcal{A})}^2\big) \\ + 6\lambda\tau(q_{k-1}f_{k-1}q_{k-1} - q_{k}f_{k}q_{k}).
\end{multline*}
Therefore, by orthogonality of martingale differences and the previous estimate, summation over $k$ gives
\begin{align*}
I  &= \sum_{k \geq 1} \| q_kf_kq_k - \mathsf{E}_{k-1}(p_{k}f_kp_{k})\|_{L_2(\mathcal{A})}^2 
\\
&\leq \lim_{k \to \infty} \left(2\left(\|q_kf_kq_k\|_{L_2(\mathcal{A})}^2 - \|q_{0}f_{0}q_{0}\|_{L_2(\mathcal{A})}^2\right) + 6\lambda \tau(q_{0}f_{0}q_{0} - q_kf_kq_k)\right) 
\\
&\leq  \lim_{k \to \infty}\left(2\|q_kf_kq_k\|_{L_2(\mathcal{A})}^2+ 6\lambda \tau(q_0f_0) \right)
\leq 8\lambda\|f\|_{L_1(\mathcal{A})},
\end{align*}
where H\"older's inequality and \eqref{item:Cuest} of Cuculescu's construction were used. To estimate $II$ we perform the telescopic sum in order to get
\begin{align*}
II \leq 2\|qfq\|_{L_2(\mathcal{A})}^2 + 2\|q_0f_0q_0\|_{L_2(\mathcal{A})}^2 \leq 4\lambda\|f\|_{L_1(\mathcal{A})},
\end{align*}
which follows from the estimate $qfq \leq \lambda q$, which in turn can be deduced from Cuculescu's construction (see \cite[Section 4.1]{parcet2009}). By this last estimate and using that $(a + b + c)^2 \leq 3a^2 + 3b^2 + 3c^2$ for $a,b,c$ positive numbers, we finally obtain
\begin{equation*}
\|g_{\Delta}\|_{L_2(\mathcal{A})}^2 \leq 39\lambda\|f\|_{L_1(\mathcal{A})}. 
\end{equation*}

The bad terms are easier to handle. Clearly the bad term $b_{\Delta}$ is comprised of the self-adjoint terms $b_k = p_{k}(f - f_k)p_{k}$ with the mean zero property $\mathsf{E}_k(b_k) = 0$, so that $\int_{\R^d} b \,d\mu = 0$. Moreover, by the orthogonality of the projections $p_k$, the tracial property of $\tau$ and since conditional expectations are bimodular and trace preserving, we have that
\begin{multline*}
\|b_{\Delta}\|_{L_1(\mathcal{A})} = \sum_{k \geq 1} \|b_k\|_{L_1(\mathcal{A})} \leq \sum_{k \geq 1}\tau\bigl(p_k(f + f_k)p_k\bigr) 
\\
= 2\tau(f(1_\mathcal{A} - q)) \leq 2\|f\|_{L_1(\mathcal{A})}.
\end{multline*}

Similarly, $\beta_{\Delta} = \sum_k\beta_k$, where $\beta_k = \mathsf{D}_{k}(p_{k} f_k p_{k}) = \mathsf{D}_k\beta_{\Delta}$ is a $k$-th martingale difference --- and hence of mean zero. Moreover, as conditional expectations are contractive on $L_1(\mathcal{A})$
\begin{multline*}
\|\beta_{\Delta}\|_{L_1(\mathcal{A})} \leq \sum_{k \geq 1} \|\beta_k\|_{L_1(\mathcal{A})} \leq 2\sum_{k \geq 1} \tau(p_{k} f_k p_{k}) = 2\tau(f(1_\mathcal{A} - q)) \leq 2\|f\|_{L_1(\mathcal{A})}.
\end{multline*}


We now turn to the off-diagonal terms, which require some more work. To get the appropriate estimate for $g_{\mathrm{off}}$, first we need to obtain a manageable expression for its $k$-th martingale difference. Rewrite $g_{\mathrm{off}}$ as
\begin{equation*}
g_{\mathrm{off}} = (1_\mathcal{A} - q)fq + qf(1_\mathcal{A} - q) + \sum_{k \geq 1}\sum_{h \geq 1} \mathsf{E}_{k+h-1}(p_kf_{k+h}p_{k+h} + p_{k+h}f_{k+h}p_k).
\end{equation*}
Since $p_{i\wedge j}, p_{i\vee j} \leq q_{i\wedge j-1}$ and by the commutation property \eqref{item:Cucomm} of  Cuculescu's construction we have that 
\begin{equation}
\label{eq:nullav}
p_if_{i\wedge j}p_j = p_iq_{i\wedge j-1}f_{i\wedge j}q_{i\wedge j-1}p_j = 0, \quad i\neq j, \quad i,j \in \mathbb{N} \cup \{\infty\}.
\end{equation}
Thus,
\begingroup
\allowdisplaybreaks
\begin{align*}
\sum_{k \geq 1}\sum_{h \geq 1} \mathsf{E}_{k+h-1}&(p_kf_{k+h}p_{k+h} + p_{k+h}f_{k+h}p_k)
\\
&= \sum_{k \geq 1}\sum_{h \geq 1} \mathsf{E}_{k+h-1}\bigl(p_k(f_{k+h} - f_k)p_{k+h} + p_{k+h}(f_{k+h} - f_k)p_k\bigr)
\\
&= \sum_{k \geq 1}\sum_{h \geq 1}\sum_{i=1}^h \mathsf{E}_{k+h-1}(p_kdf_{k+i}p_{k+h} + p_{k+h}df_{k+i}p_k).
\end{align*}
We may now proceed to calculate $\mathsf{D}_j(g_{\mathrm{off}})$ for $j \geq 1$. Taking into account that, for $h \geq 1$, $\mathsf{D}_j\mathsf{E}_{k+h-1} = \mathsf{D}_j$ if $j < k + h$ and zero otherwise, we get that
\begin{multline*}
\mathsf{D}_j(g_{\mathrm{off}}) = \mathsf{D}_j((1_\mathcal{A} - q)fq + qf(1_\mathcal{A} - q))
\\
+ \sum_{k<j}\,\sum_{h > j-k} \, \sum_{i=1}^h \mathsf{D}_j(p_kdf_{k+i}p_{k+h} + p_{k+h}df_{k+i}p_k)
\\
+ \sum_{k\geq j}\,\sum_{h \geq 1} \, \sum_{i=1}^h \mathsf{D}_j(p_kdf_{k+i}p_{k+h} + p_{k+h}df_{k+i}p_k) = I + II + III.
\end{multline*}
We deal first with $II$. By Fubini's theorem we obtain that
\begingroup
\allowdisplaybreaks
\begin{align*}
II &= \sum_{k<j} \Bigg(\sum_{i=1}^{j-k} \, \sum_{h > j-k} \mathsf{D}_j(p_kdf_{k+i}p_{k+h} + p_{k+h}df_{k+i}p_k) 
\\
& \phantom{\sum_{k<j} \Bigg(\sum_{i=1}^{j-k} \, \sum_{h > j-k} \mathsf{D}_j}
+ \sum_{i > j-k} \, \sum_{h \geq i}   \mathsf{D}_j(p_kdf_{k+i}p_{k+h} + p_{k+h}df_{k+i}p_k)\Bigg)
\\
&= \sum_{k<j} \Bigg(\sum_{i=1}^{j-k} \mathsf{D}_j(p_kdf_{k+i}q_j + q_jdf_{k+i}p_k) 
\\
&\phantom{\sum_{k<j} \Bigg(\sum_{i=1}^{j-k} \mathsf{D}_j}
+ \sum_{i > j-k} \mathsf{D}_j(p_kdf_{k+i}q_{k+i-1} + q_{k+i-1}df_{k+i}p_k)
\\
&\phantom{\phantom{\sum_{k<j} \Bigg(\sum_{i=1}^{j-k} \mathsf{D}_j} + \sum_{i > j-k}}
- \sum_{i \geq 1} \mathsf{D}_j(p_kdf_{k+i}q + qdf_{k+i}p_k)\Bigg) = II_1 + II_2 + II_3.
\end{align*}
After summing over $i$ in $II_1$ and noticing that by \eqref{item:Cucomm} of Cuculescu's construction (recall that $k < j$)
\begin{equation*}
p_kf_kq_j = p_kq_{k-1}f_kq_{k-1}q_j = 0 = q_jf_kp_k,
\end{equation*}
we find that
\begin{equation*}
II_1 = \sum_{k < j} \mathsf{D}_j(p_kf_jq_j + q_jf_jp_k) = \mathsf{D}_j\bigl((1_\mathcal{A} - q_{j-1})f_jq_j + q_jf_j(1_\mathcal{A} - q_{j-1})\bigr).
\end{equation*}
The term $II_2$ vanishes since
\begin{equation}
\label{eq:Dmvanish}
p_kdf_{k + i}q_{k + i-1} + q_{k+i-1}df_{k+i}p_k = \mathsf{D}_{k+i}(p_kfq_{k+i-1} + q_{k+i-1}fp_k)
\end{equation}
and $\mathsf{D}_j\mathsf{D}_{k+i} = 0$, as $k + i > j$. Performing the summation over $i$ in $II_3$ and using \eqref{eq:nullav} with $i \wedge j = k$ and $i \vee j = \infty$, we get that
\begin{multline*}
II = \mathsf{D}_j\bigl((1_\mathcal{A} - q_{j-1})f_jq_j + q_jf_j(1_\mathcal{A} - q_{j-1})\bigr)
\\ 
- \mathsf{D}_j\bigl((1_\mathcal{A} - q_{j-1})fq + qf(1_\mathcal{A} - q_{j-1})\bigr).
\end{multline*}
Changing the order of summation
\begingroup
\allowdisplaybreaks
\begin{align*}
III &= \sum_{k \geq j}\,\sum_{i \geq 1} \, \sum_{h\geq i} \mathsf{D}_j(p_kdf_{k+i}p_{k+h} + p_{k+h}df_{k+i}p_k)
\\
&= \sum_{k \geq j} \Bigg(\,\sum_{i \geq 1} \mathsf{D}_j(p_kdf_{k+i}q_{k+i-1} + q_{k+i-1}df_{k+i}p_k)
\\
&\phantom{= \sum_{j\geq k} \Bigg(\,\sum_{i \geq 1}}
- \sum_{i \geq 1} \mathsf{D}_j(p_kdf_{k+i}q + qdf_{k+i}p_k) \Bigg)
\\
&= \minus\mathsf{D}_j\bigl((q_{j-1} - q)fq + qf(q_{j-1} - q)\bigr).
\end{align*}
Here, we have also used \eqref{eq:Dmvanish}, as $k+i > j$, and \eqref{eq:nullav} with $i \vee j = \infty$. Finally, summing everything we get that for $j \geq 1$
\begin{equation*}
\label{eq:dkgoff}
\mathsf{D}_j(g_{\mathrm{off}}) = \mathsf{D}_j\bigl((1_{\mathcal{A}} - q_{j-1})f_jq_j) + \mathsf{D}_j(q_jf_j(1_{\mathcal{A}} - q_{j-1})\bigr).
\end{equation*}
On the other hand, $\mathsf{D}_j(g_{\mathrm{off}})=0$ for $j \leq 0$. Indeed,
\begin{multline*}
\mathsf{D}_j(g_{\mathrm{off}}) = \mathsf{D}_j((1_\mathcal{A} - q)fq + qf(1_\mathcal{A} - q))
\\
+ \sum_{k\geq 1}\,\sum_{h \geq 1} \, \sum_{i=1}^h \mathsf{D}_j(p_kdf_{k+i}p_{k+h} + p_{k+h}df_{k+i}p_k)
\end{multline*}
and, arguing as with $III$ above and since $q_0 = 1_\mathcal{A}$, we have that
\begin{equation*}
 \sum_{k\geq 1}\,\sum_{h \geq 1} \, \sum_{i=1}^h \mathsf{D}_j(p_kdf_{k+i}p_{k+h} + p_{k+h}df_{k+i}p_k) = \minus\mathsf{D}_j\bigl((1_\mathcal{A} - q)fq + qf(1_\mathcal{A} - q)\bigr).
\end{equation*}
Thus, in $L_2$ sense
\begingroup
\allowdisplaybreaks
\begin{multline*}
g_{\mathrm{off}} = \sum_{j\geq 1} \mathsf{D}_j (g_{\mathrm{off}})
= \sum_{j \geq 1} \sum_{k < j}\mathsf{D}_j(p_kf_jq_j + q_jf_jp_k) \notag
\\
= \sum_{k \geq 1}\sum_{h \geq 1} \mathsf{D}_{k + h}(p_kf_{k+h}q_{k+h} + q_{k+h}f_{k+h}p_k)
=: \sum_{k \geq 1}\sum_{h \geq 1} g_{k,h}.
\end{multline*}
We are now in the position to prove the estimate in \eqref{item:LdosCZoff} of Theorem \ref{thm:CZd}. Notice first that by H\"older's inequality, the $C^*$-algebra property and \eqref{item:Cuest} of Cuculescu's construction
\begingroup
\allowdisplaybreaks
\begin{align*}
\|g_{k,h}\|_{L_2(\mathcal{A})}^2 &\leq 16 \| q_{k+h}f_{k+h}p_k \|_{L_2(\mathcal{A})}^2 
\\
&= 16\tau(p_kf_{k+h}q_{k+h}f_{k+h}p_k) 
\\
&\leq 16 \big\|f_{k+h}^{1/2} q_{k+h}f_{k+h}^{1/2} \big\|_\mathcal{A} \,\tau\bigl(f_{k+h}^{1/2} p_k f_{k+h}^{1/2}\bigr) 
\\
&= 16 \|q_{k+h} f_{k+h} q_{k+h} \|_\mathcal{A} \, \tau(p_kf_{k+h}p_k) 
\leq 16\lambda \tau(fp_k).
\end{align*}
This proves that for all $h\geq 1$
\begin{equation*}
\sum_{k \geq 1} \|g_{k,h}\|_{L_2(\mathcal{A})}^2 \leq 16\lambda \tau(f(1_\mathcal{A} - q)) \leq  16\lambda\|f\|_{L_1(\mathcal{A})}.
\end{equation*}
For the bad terms we follow \cite{parcet2009}. First, rewrite $b_{\mathrm{off}}$ as 
\begin{equation*}
b_{\mathrm{off}} = \sum_{h \geq 1}\sum_{k \geq 1}p_k(f - f_{k+h})p_{k+h} + p_{k+h}(f - f_{k+h})p_k =: \sum_{h \geq 1}\sum_{k \geq 1}b_{k,h}.
\end{equation*}
Clearly, the terms $b_{k,h}$ have mean zero and satisfy the estimate
\begin{equation*}
\|b_{k,h}\|_{L_1(\mathcal{A})} \leq 2 \|p_k f p_{k+h} + p_{k+h} f p_k \|_{L_1(\mathcal{A})}.
\end{equation*}
Next, observe that we can decompose the off-diagonal terms $p_k f p_{k+h} + p_{k+h} f p_k$ into a sum of four positive overlapping box-diagonal terms
\begingroup
\allowdisplaybreaks
\begin{multline*}
p_k f p_{k+h} + p_{k+h} f p_k  = \Biggl(\sum_{j = 0}^h p_{k + j}\Biggr)f\Biggl(\sum_{j = 0}^h p_{k + j}\Biggr)
- \Biggl(\sum_{j = 0}^{h-1} p_{k + j}\Biggr)f\Biggl(\sum_{j = 0}^{h-1} p_{k + j}\Biggr)
\\
- \Biggl(\sum_{j = 1}^h p_{k + j}\Biggr)f\Biggl(\sum_{j = 1}^h p_{k + j}\Biggr)
+ \Biggl(\sum_{j = 1}^{h-1} p_{k + j}\Biggr)f\Biggl(\sum_{j = 1}^{h-1} p_{k + j}\Biggr).
\end{multline*}
The previous expression implies that
\begingroup
\allowdisplaybreaks
\begin{align*}
\sum_{k \geq 1} \|p_k f p_{k+h} + p_{k+h} f p_k \|_{L_1(\mathcal{A})} 
&\leq 4 \sum_{k \geq 1} \sum_{j=0}^h\tau(fp_{k+j})
\\
&= 4 \sum_{j=0}^h  \tau(f(q_j - q)) \leq 4(h+1)\|f\|_{L_1(\mathcal{A})},
\end{align*}
and hence the estimate in \eqref{item:badbCZoff} holds. On the other hand, we have
\begin{equation*}
\beta_{\mathrm{off}} = \sum_{k \geq 1}\sum_{h\geq 1} \mathsf{D}_{k+h}(p_{k}f_{k+h}p_{k+h} + p_{k+h}f_{k+h}p_{k}) =: \sum_{k \geq 1}\sum_{h\geq 1} \beta_{k,h}.
\end{equation*}
Each term in the previous sum satisfies the same estimate
\begin{equation*}
\|\beta_{k,h}\|_{L_1(\mathcal{A})} \leq 2 \|p_k f p_{k+h} + p_{k+h} f p_k \|_{L_1(\mathcal{A})},
\end{equation*}
which yields the corresponding estimate for $\beta_{\mathrm{off}}$. \qed


\section{Commuting Haar shift operators}
\label{section:CHS}

We now turn to the proof of Theorem \ref{thm:HS}. Namely that
\begin{equation*}
\lambda \tau\left( \left\{\left|\Sha_{r,s} f \right| > \lambda \right\} \right) \lesssim \|f\|_{L_1(\mathcal{A})}
\end{equation*}
for all $\lambda > 0$. Here $\tau (\{ |f| > \lambda \})$ denotes the trace of the spectral projection of $|f|$ associated to the interval $(\lambda,\infty)$, which defines a noncommutative distribution function. We find this terminology more intuitive, since it is reminiscent of the classical one. Following the construction of noncommutative symmetric spaces (see \cite{mei-parcet2009} and references therein), the resulting $L_{1,\infty}(\mathcal{A})$ space is a quasi-Banach space with quasi-norm $\|f\|_{L_{1,\infty}(\mathcal{A})} = \sup_{\lambda > 0} \lambda \tau (\{ |f| > \lambda \})$ which interpolates with $L_2(\mathcal{A})$. It should be mentioned that the weak Bochner space $L_{1,\infty}(\R^d, \mu; L_1(\mathcal{M}))$ is of no use for our purposes since $L_1(\mathcal{M})$ is not a $\mathrm{UMD}$ space and thus even Haar multipliers may not be bounded, which rules out the use of this space as an appropriate setting for  providing weak-type $(1,1)$ estimates for the operators in question. The same applies if one considers $\mathcal{M}$ instead of $L_1(\mathcal{M})$ as target space.

\subsection*{Proof of Theorem \ref{thm:HS}}
Let $f \in \mathcal{A}_{+,K}$. The general case follows by the density of the span of $\mathcal{A}_{+,K}$ in $L_1(\mathcal{A})$. Consider the Calder\'on-Zygmund decomposition $f = g_{\Delta} + b_{\Delta} + \beta_{\Delta} + g_{\mathrm{off}}  + b_{\mathrm{off}} + \beta_{\mathrm{off}}$ associated to $(f, \lambda)$ for a given $\lambda > 0$. By the quasi-triangle inequality in $L_{1,\infty}(\mathcal{A})$ it suffices to show that
\begin{equation*}
\lambda \tau (\{|\Sha_{r,s}(\gamma)| > \lambda\}) \lesssim \|f\|_{L_1(\mathcal{A})}
\end{equation*}
for all $\gamma \in \{g_{\Delta}, b_{\Delta}, \beta_{\Delta}, g_{\mathrm{off}}, b_{\mathrm{off}}, \beta_{\mathrm{off}}\}$. We start with the diagonal terms, for which the estimates are very similar to the classical ones. For $g_{\Delta}$ we use Chebyshev's inequality, the $L_2$ boundedness of $\Sha_{r,s}$ and the $L_2$ estimate in \eqref{item:LdosCZd} of Theorem \ref{thm:CZd} to get
\begin{equation*}
\lambda \tau \{|\Sha_{r,s}(g_{\Delta})| > \lambda\} 
\leq 39\|\Sha_{r,s}\|_{\mathcal{B}(L_2(\mathcal{A}))}^2 \|f\|_{L_1(\mathcal{A})},
\end{equation*} 
where $\|\Sha_{r,s}\|_{\mathcal{B}(L_2(\mathcal{A}))}$ denotes the operator norm of $\Sha_{r,s}$ on $L_2(\mathcal{A})$. For the remaining $\gamma$, we decompose $\Sha_{r,s}(\gamma)$ as 
\begin{multline*}
\Sha_{r,s}(\gamma) = (1_\mathcal{A}-q)\Sha_{r,s}(\gamma)(1_\mathcal{A} - q) +  q\Sha_{r,s}(\gamma)q  
\\
+ q\Sha_{r,s}(\gamma)(1_\mathcal{A} - q) + (1_\mathcal{A}-q)\Sha_{r,s}(\gamma)q.
\end{multline*} 
Since the distribution function is adjoint-invariant and by the second estimate in \eqref{item:Cucuest} of Cuculescu's construction, we get that
\begin{align*}
\lambda \tau (\{|\Sha_{r,s}(\gamma)| > \lambda\}) &\leq 12\|f\|_{L_1(\mathcal{A})} + \lambda\tau(\{|q\Sha_{r,s}(\gamma)q| > \lambda/4\}).
\end{align*}
To prove the estimate for $\gamma = b_{\Delta}$, observe that we may further decompose each term $b_k$ in \eqref{item:badbCZd} of Theorem \ref{thm:CZd} as 
\begin{align*}
b_k = \sum_{L \in \mathscr{D}_k} p_L(f - \langle f \rangle_L)p_L 1_L =: \sum_{L \in \mathscr{D}_k} b_L,
\end{align*}
where the projections $p_L$ are defined as in \eqref{eq:CZcubes}. Since the Haar function $\phi_R$ is constant on dyadic subcubes of $R$ and $b_L$ has zero integral, $\langle b_L, \phi_R \rangle$ is nonzero only for $R \subset L$, i.e., $R^{(r)} \subset L^{(r)}$ for their respective $r$-dyadic ancestors. On the other hand, if $x\in L$ we have that $q(x) \leq q_k(x) = q_L$  in the order of the lattice $\mathcal{P}(\mathcal{M})$. This together with \eqref{eq:CZcubes} gives that for $x\in L$
\begin{equation}
\label{eq:suppShabL}
q(x) \langle b_L, \phi_R \rangle q(x) = q(x) q_L \, p_L\langle b_L, \phi_R \rangle p_L q_L q(x) = 0.
\end{equation}
Using that  $\alpha^Q_{R,S} \in \mathcal{M} \cap \mathcal{M}'$ we find the estimate
\begingroup
\allowdisplaybreaks
\begin{align}
\| q &\Sha_{r,s} (b_L) q \|_{L_1(\mathcal{A})} \label{eq:Shabdest}
\\*
&\leq \sum_{\begin{subarray}{c} Q \in \mathscr{D} \\ L \subsetneq Q \subset L^{(r)} \end{subarray}} \sum_{\begin{subarray}{c} R \in \mathscr{D} _r(Q), \, R \subset L\\ S \in \mathscr{D} _s(Q) \end{subarray}} 
\|\alpha_{R,S}^Q\|_\mathcal{M} \|\langle b_L, \phi_R \rangle\|_{L_1(\mathcal{M})} \|\psi_S\|_{L_1(\mu)} \notag
\\ 
&\leq \sup_{Q, R, S} \|\alpha_{R,S}^Q\|_\mathcal{M} \sum_{\begin{subarray}{c} Q \in \mathscr{D} \\ L \subsetneq Q \subset L^{(r)} \end{subarray}} \sum_{\begin{subarray}{c} R \in \mathscr{D} _r(Q), \, R \subset L\\ S \in \mathscr{D} _s(Q) \end{subarray}} \| \phi_R \|_{L_{\infty}(\mu)} \| \psi_S \|_{L_1(\mu)} \| b_L\|_{L_1(\mathcal{A})} \notag
\\
&\leq r 2^{(r+s)d}\sup_{Q, R, S} \|\alpha_{R,S}^Q\|_\mathcal{M}  \,\, \Xi(\Phi, \Psi; r,s) \| b_L\|_{L_1(\mathcal{A})} \notag.
\end{align}
This, Chebyshev's inequality, the fact that dyadic cubes in $\mathscr{D}_k$ are disjoint and \eqref{item:badbCZd} of Theorem \ref{thm:CZd} give the estimate
\begin{equation*}
\lambda \tau ( \{ |q\Sha_{r,s} b_{\Delta} q| > \lambda \}) \leq r 2^{1 + (r+s)d}\sup_{Q, R, S} \|\alpha_{R,S}^Q\|_\mathcal{M}  \;  \Xi(\Phi, \Psi; r,s) \|  f \|_{L_1(\mathcal{A})}.
\end{equation*}
For $\gamma = \beta_{\Delta}$ we proceed likewise by writing
\begin{align*}
\beta_k &= \mathsf{D}_k(\beta_{\Delta}) = \sum_{L \in \mathscr{D}_{k-1}} \sum_{J \in \mathscr{D}_1(L)} p_J \langle f \rangle_J \, p_J \Biggl(1_J - \frac{\mu(J)}{\mu(L)}1_{L}  \Biggr) 
\\*
&=: \sum_{L \in \mathscr{D}_{k-1}}\sum_{J \in \mathscr{D}_1(L)}  \beta_{L,J}=: \sum_{L \in \mathscr{D}_k} \beta_L,
\end{align*}
where each term $\beta_L$ is supported (as an operator-valued function) on $L$, is constant on the dyadic descendants of $L$ and has mean zero. By  Chebyshev's inequality we have
\begin{multline*}
\lambda\tau(\{|q\Sha_{r,s}(\beta_{\Delta})q| > \lambda\})
\leq 
\sum_{k \geq 1}\sum_{L \in \mathscr{D}_{k-1}} \Bigl( \int_{\R^d\setminus L} \nu \bigl(|q(x)\Sha_{r,s}  \beta_L(x) q(x)|\bigr) d\mu(x) 
\\
+ \int_{L} \nu \bigl(|q(x)\Sha_{r,s}  \beta_L(x)q(x)|\bigr) d\mu(x)\Bigr).
\end{multline*}
Since $\langle \beta_L, \phi_R \rangle$ is nonzero only for dyadic cubes $R \subset L$, proceeding as in \eqref{eq:Shabdest} we obtain
\begin{multline*}
\int_{\R^d\setminus L} \nu \bigl(|q(x)\Sha_{r,s}  \beta_L(x) q(x)|\bigr) d\mu(x) 
\\*
\leq r 2^{(r+s)d}\sup_{Q, R, S} \|\alpha_{R,S}^Q\|_\mathcal{M}  \;  \Xi(\Phi, \Psi; r,s) \|  \beta_L \|_{L_1(\mathcal{A})}.
\end{multline*}
Arguing as above and recalling that $\alpha^Q_{R,S} \in \mathcal{M} \cap \mathcal{M}'$, for $x \in L$ we obtain
\begingroup
\allowdisplaybreaks
\begin{align*}
q(x) \Sha_{r,s} (\beta_L)(x) q(x)
& = 
\sum_{\begin{subarray}{c} Q \in \mathscr{D} \\ L \subset Q \subset L^{(r)} \end{subarray}} \; \sum_{\begin{subarray}{c} R \in \mathscr{D}_r(Q), \, R \subset L\\ S \in \mathscr{D} _s(Q) \end{subarray}} 
\alpha_{R,S}^Q \, q(x) \langle \beta_L, \phi_R \rangle q(x) \psi_S(x) 
\\ 
& \phantom{\sum_{\begin{subarray}{c} Q \in \mathscr{D} \\ L \subset Q \subset L^{(r+1)} \end{subarray}} \;}  
+ \sum_{\begin{subarray}{c} Q \in \mathscr{D} \\ Q \subsetneq L \end{subarray}} \; \sum_{\begin{subarray}{c} R \in \mathscr{D}_r(Q)\\ S \in \mathscr{D} _s(Q) \end{subarray}} 
\alpha_{R,S}^Q \, q(x) \langle \beta_L, \phi_R \rangle q(x) \psi_S(x)
\\ 
&= F_L(x) + G_L(x).
\end{align*}
As in \eqref{eq:Shabdest} we get the estimate
\begin{equation*}
\int_{L} \nu \bigl(|F_L(x)|\bigr) d\mu(x) 
\leq (r+1) 2^{(r+s)d}\sup_{Q, R, S} \|\alpha_{R,S}^Q\|_\mathcal{M}  \;  \Xi(\Phi, \Psi; r,s) \|  \beta_L \|_{L_1(\mathcal{A})}.
\end{equation*}
To estimate $G_L(x)$ we further decompose $\beta_L$ and get
\begin{equation*}
G_L(x) = \sum_{J \in \mathscr{D}_1(L)} \sum_{\begin{subarray}{c} Q \in \mathscr{D} \\ Q \subsetneq L \end{subarray}} \; \sum_{\begin{subarray}{c} R \in \mathscr{D}_r(Q)\\ S \in \mathscr{D} _s(Q) \end{subarray}} 
\alpha_{R,S}^Q \, q(x) \langle \beta_{L,J}, \phi_R \rangle q(x) \psi_S(x)
= \sum_{J \in \mathscr{D}_1(L)} G_{L,J}(x).
\end{equation*}
Given $J \in \mathscr{D}_1(L)$ and a dyadic cube $Q \subsetneq L$ we either have $Q \subset J$ or $Q \subset L \setminus J$. Yet the former case  leads to zero terms since, as in \eqref{eq:suppShabL}, for $x \in Q \subset J$ we have $q(x) \leq  q_J$ and thus $q(x) \langle \beta_{L,J}, \phi_R \rangle q(x) =  0$. Hence,
\begingroup
\allowdisplaybreaks
\begin{align*}
G_{L,J}&(x) = \minus p_J \langle f \rangle_J \,p_J \, \frac{\mu(J)}{\mu(L)} \sum_{\begin{subarray}{c} Q \in \mathscr{D} \\ Q \subset L \setminus J \end{subarray}} \; \sum_{\begin{subarray}{c} R \in \mathscr{D}_r(Q)\\ S \in \mathscr{D} _s(Q) \end{subarray}} 
\alpha_{R,S}^Q \, q(x) \langle 1_{L \setminus J}, \phi_R \rangle q(x) \psi_S(x)
\\
&= \minus p_J \langle f \rangle_J \,p_J \, \frac{\mu(J)}{\mu(L)} \sum_{\begin{subarray}{c} Q' \in \mathscr{D}_1(L) \\ Q' \neq J \end{subarray}} \;\sum_{Q \in \mathscr{D}(Q')} \; \sum_{\begin{subarray}{c} R \in \mathscr{D}_r(Q)\\ S \in \mathscr{D} _s(Q) \end{subarray}} 
\alpha_{R,S}^Q \, q(x) \langle 1_{Q'}, \phi_R \rangle q(x) \psi_S(x)
\\
&= \minus p_J \langle f \rangle_J \,p_J \, \frac{\mu(J)}{\mu(L)}\sum_{\begin{subarray}{c} Q' \in \mathscr{D}_1(L) \\ Q' \neq J \end{subarray}} q(x)\Sha_{r,s}^{Q'}(1_{Q'})(x)q(x).
\end{align*}
Then, by H\"older's inequality and the fact that $\mathrm{supp}_{\R^d} \bigl(\Sha_{r,s}^{Q'}(1_{Q'})\bigr) \subset Q'$
\begingroup
\allowdisplaybreaks
\begin{align*}
\int_{L} & \nu \bigl(|G_L(x)|\bigr) d\mu(x)
\\*
&=  \int_{L} \nu \Biggl(\,\Bigg|\sum_{J \in \mathscr{D}_1(L)} p_J \langle f \rangle_J \,p_J \, \frac{\mu(J)}{\mu(L)} \,\sum_{\begin{subarray}{c} Q' \in \mathscr{D}_1(L) \\ Q' \neq J \end{subarray}} q(x)\Sha_{r,s}^{Q'}(1_{Q'})(x)q(x) \Bigg|\,\Biggr)d\mu(x)
\\
&\leq \sum_{J \in \mathscr{D}_1(L)} \|p_J \langle f \rangle_J \, p_J\|_{L_1(\mathcal{M})} \, \frac{\mu(J)}{\mu(L)} \sum_{\begin{subarray}{c} Q' \in \mathscr{D}_1(L) \\ Q' \neq J \end{subarray}} \int_{L} \|\Sha_{r,s}^{Q'}(1_{Q'})(x)\|_{\mathcal{M}}\,d\mu(x)
\\
&\leq \sum_{J \in \mathscr{D}_1(L)} \|p_J \langle f \rangle_J \, p_J\|_{L_1(\mathcal{M})} \, \frac{\mu(J)}{\mu(L)}
\\
& \qquad \qquad \qquad \qquad \qquad \times \left(\sum_{\begin{subarray}{c} Q' \in \mathscr{D}_1(L) \\ Q' \neq J \end{subarray}} \left( \int_{\R^d} \|\Sha_{r,s}^{Q'}(1_{Q'})(x)\|_{\mathcal{M}}^2\,d\mu(x)\right)^{\frac{1}{2}} \mu(Q')^{\frac{1}{2}}\right)
\\
&\leq \sup_{\begin{subarray}{c} Q \in \mathscr{D}, \\ \mu(Q) \neq 0 \end{subarray}} \frac{1}{\mu(Q)^{\frac{1}{2}}}\left( \int_{\R^d} \|\Sha_{r,s}^{Q}(1_{Q})(x)\|_{\mathcal{M}}^2\,d\mu(x)\right)^{\frac{1}{2}} \, \Bigg\| \sum_{J \in \mathscr{D}_1(L)} p_J \langle f \rangle_J \, p_J 1_J\Bigg\|_{L_1(\mathcal{A})},
\end{align*}
which is finite by the local vector-valued $L_2$ estimate \eqref{eq:localv-v}. By the estimate in \eqref{item:badbetaCZd} of the Calder\'{o}n-Zygmund decomposition
\begin{multline*}
\sum_{k \geq 1}\sum_{L \in \mathscr{D}_k} \Biggl( \| \beta_L \|_{L_1(\mathcal{A})} + \Bigg\| \sum_{J \in \mathscr{D}_1(L)} p_J \langle f \rangle_J \, p_J 1_J\Bigg\|_{L_1(\mathcal{A})} \, \Biggr) 
\\[5pt]
 \leq \sum_k \bigl(\|\beta_k\|_1 + \|p_k f_k p_k\|_1\bigr)  \leq 3 \|f\|_1.
\end{multline*}
Thus, gathering the previous estimates
\begin{align*}
\lambda\tau(\{|q\Sha_{r,s}&(\beta_{\Delta})q| > \lambda\})
\\*
&\leq  \Biggl( (r+2) 2^{1+(r+s)d}\sup_{Q, R, S} \|\alpha_{R,S}^Q\|_\mathcal{M}  \;  \Xi(\Phi, \Psi; r,s) 
\\
& \quad\quad + \sup_{\begin{subarray}{c} Q, \in \mathscr{D} \\ \mu(Q) \neq 0 \end{subarray}} \frac{1}{\mu(Q)^{\frac{1}{2}}}\left( \int_{\R^d} \|\Sha_{r,s}^{Q}(1_{Q})(x)\|_{\mathcal{M}}^2\,d\mu(x)\right)^{\frac{1}{2}} \Biggr)\|f\|_{L_1(\mathcal{A})}.
\end{align*}


We now turn to the weak-type estimates for the off-diagonal terms, starting with $g_{\mathrm{off}}$. By Chebyshev's inequality 
\begin{equation*}
\lambda \tau ( \{ |q \Sha_{r,s} (g_{\mathrm{off}}) q| > \lambda \}) \leq \frac{1}{\lambda}\left(\, \sum_{h \geq 1} 
\Bigg\| \sum_{k \geq 1} q \Sha_{r,s}(g_{k,h}) q \Bigg\|_{L_2(\mathcal{A})} \right)^2.
\end{equation*}
We further decompose the terms $g_{k,h}$ as
\begin{equation*}
g_{k,h} = \sum_{L \in \mathscr{D}_k}\sum_{J \in  \mathscr{D}_h(L)} \bigl(p_L\langle f \rangle_J\,q_J + q_J\langle f \rangle_J\,p_L \bigr) \Biggl(1_J - \frac{\mu(J)}{\mu(\widehat{J})}1_{\widehat{J}}  \Biggr) =: \sum_{L \in \mathscr{D}_k} g_{L,h}.
\end{equation*}
Clearly, each term $g_{L,h}$ is such that $\mathrm{supp}_{\R^d}(g_{L,h}) \subset L$ and has mean zero on the $(h-1)$-descendants of $L$. Thus, $\langle g_{L,h} , \phi_R \rangle$ is nonzero only for $R  \subset \widehat{J}$ for some $J \in \mathscr{D}_h(L)$, which amounts to say that $R \in \mathscr{D}_{h+j-1}(L)$ for some $j \geq 0$. Furthermore  $g_{L,h} = p_LA_{L,h} + A_{L,h}^*p_L$, where
\begin{equation*}
 A_{L,h} = p_L\langle f \rangle_J\,q_J \Biggl(1_J - \frac{\mu(J)}{\mu(\widehat{J})}1_{\widehat{J}}  \Biggr).
\end{equation*}
Proceeding as in \eqref{eq:suppShabL} we get that $q(x) \langle g_{L,h}, \phi_R \rangle q(x) = 0$ if $x \in L$. In other words, only the cubes $R$ such that $R^{(r)} \supsetneq L$ lead to nonzero terms. These two observations in terms of side-lengths provide that $h$ must be such that $\ell(L) = 2^{-k} < \ell(R^{(r)}) = 2^{-(k+h+j-1-r)}$, namely $h \leq r$. This and the assumption $\alpha^Q_{R,S} \in \mathcal{M} \cap \mathcal{M}'$ allow us to deduce that $q(x)\Sha_{r,s}g_{k,h} (x) q(x)=0$ whenever $h > r$. This localization property and the orthogonality of martingale differences in $L_2(\mathcal{A})$, enable us to obtain that
\begin{align*}
\sum_{h \geq 1} 
\Bigg\| \sum_{k \geq 1} q \Sha_{r,s}(g_{k,h}) q \Bigg\|_{L_2(\mathcal{A})}
&\leq \|\Sha_{r,s}\|_{\mathcal{B}(L_2(\mathcal{A}))} \sum_{h=1}^r  \Bigg\|\sum_{k \geq 1} g_{k,h} \Bigg\|_{L_2(\mathcal{A})}
\\
&= \|\Sha_{r,s}\|_{\mathcal{B}(L_2(\mathcal{A}))} \sum_{h=1}^r \Biggl(\,\sum_{k \geq 1} \|g_{k,h}\|_{L_2(\mathcal{A})}^2\Biggr)^{1/2} .
\end{align*}
Therefore, by the estimate in \eqref{item:LdosCZoff} of Theorem \ref{thm:CZd} we arrive at
\begin{equation*}
\lambda \tau ( \{ |q \Sha_{r,s} (g_{\mathrm{off}}) q| > \lambda \}) \leq 16 r^2 \|\Sha_{r,s}\|_{\mathcal{B}(L_2(\mathcal{A}))}^2 \|f\|_{L_1(\mathcal{A})}.
\end{equation*}

To get the estimate for $b_{\mathrm{off}}$ we proceed in an entirely similar way by decomposing the terms $b_{k,h}$ in \eqref{item:badbCZoff} of Theorem \ref{thm:CZd} as
\begin{equation*}
b_{k,h} =\sum_{L \in \mathscr{D}_k}\sum_{J \in  \mathscr{D}_h(L)} \bigl(p_L(f - \langle f \rangle_J)p_J + p_J(f - \langle f \rangle_J)p_L\bigr) 1_J =: \sum_{L \in \mathscr{D}_k} b_{L,h}.
\end{equation*}
It is clear that $\mathrm{supp}_{\R^d} (b_{L,h}) \subset L$, that $b_{L,h}$ has mean zero over the $h$-dyadic descendants of $L$ and that $b_{L,h} = p_LB_{L,h} + B_{L,h}^*p_L$, with $B_{L,h} = p_L(f - \langle f \rangle_J)p_J1_J$.  Arguing as above, $q(x)\langle b_{L,h}, \phi_R \rangle q(x)$ is nonzero only for $R  \subsetneq L \subsetneq R^{(r)} \subsetneq L^{(r)}$ and hence $q(x) \Sha_{r,s} (b_{L,h})(x) q(x)$ vanishes if  $h > r$. Thus, for $h \leq r$ we follow the steps in \eqref{eq:Shabdest} to get the estimate
\begin{align*}
\sum_{L \in \mathscr{D}_k} \| q \Sha_{r,s} (b_{L,h}) q \|_{L_1(\mathcal{A})} \leq (r-1) 2^{(r+s)d}\sup_{Q, R, S} \|\alpha_{R,S}^Q\|_\mathcal{M}  \,\, \Xi(\Phi, \Psi; r,s) \| b_{k,h}\|_{L_1(\mathcal{A})}.
\end{align*}
By Chebyshev's inequality and the estimate in \eqref{item:badbCZoff} of the Calder\'{o}n-Zygmund decomposition we obtain
\begin{align*}
\lambda \tau ( \{ |q &\Sha_{r,s} (b_{\mathrm{off}}) q| > \lambda \}) 
\\
&\leq (r-1) 2^{3+(r+s)d}\sup_{Q, R, S} \|\alpha_{R,S}^Q\|_\mathcal{M}  \,\, \Xi(\Phi, \Psi; r,s) \sum_{h=1}^{r}(h+1)\|f\|_{L_1(\mathcal{A})}
\\
&= r(r-1)(r+3) 2^{2+(r+s)d}\sup_{Q, R, S} \|\alpha_{R,S}^Q\|_\mathcal{M}  \,\, \Xi(\Phi, \Psi; r,s) \|f\|_{L_1(\mathcal{A})}.
\end{align*}
Finally, for $\gamma = \beta_{\mathrm{off}}$ observe that
\begin{align*}
\beta_{k,h} &= \sum_{L \in \mathscr{D}_k}\sum_{J \in  \mathscr{D}_h(L)} \bigl(p_L\langle f \rangle_Jp_J + p_J\langle f \rangle_Jp_L \bigr) \Biggl(1_J - \frac{\mu(J)}{\mu(\widehat{J})}1_{\widehat{J}}  \Biggr) 
\\
&=: \sum_{L \in \mathscr{D}_k} \beta_{L,h} = \sum_{L \in \mathscr{D}_k} \bigl( p_LC_{L,h} +C_{L,h}^*p_L  \bigr).
\end{align*}
Here we may repeat the analysis made for $b_{L,h}$, as each $\beta_{L,h}$ is a  $(k+h)$-martingale difference operator with $\mathrm{supp}_{\R^d}(\beta_{L,h}) \subset L$. This and \eqref{item:badbetaCZoff} of Theorem \ref{thm:CZd} render the desired estimate
\begin{align*}
\lambda \tau ( \{ |q &\Sha_{r,s} (\beta_{\mathrm{off}}) q| > \lambda \}) 
\\
&\leq r(r-1)(r+3) 2^{2+(r+s)d}\sup_{Q, R, S} \|\alpha_{R,S}^Q\|_\mathcal{M}  \,\, \Xi(\Phi, \Psi; r,s) \|f\|_{L_1(\mathcal{A})},
\end{align*}
with which we complete the proof of Theorem \ref{thm:HS}. \qed

\begin{remark}
It is worth mentioning that we have not truly needed the assumption that the symbols are commuting to obtain the estimates for the diagonal terms. Indeed, all the calculations for the diagonal terms in the proof of Theorem \ref{thm:HS} can be done without this assumption simply by rearranging multiplications. Unlike in \eqref{eq:suppShabL}, in the case when $\gamma \in \{g_{\mathrm{off}}, b_{\mathrm{off}}, \beta_{\mathrm{off}}\}$ and $x \in L$, $q(x)$ is required to be multiplied on \emph{both} sides of $\langle \gamma_{L,h}, \phi_R \rangle$ in order to annihilate it.
\end{remark}

\begin{remark}
The consideration of noncommuting symbols in \eqref{eq:ncHaar} introduces considerable additional difficulties when trying to provide a priori estimates. Firstly, different operators arise depending on whether the symbols act by right or left multiplication on each coefficient $\langle f, \phi_R \rangle$. More specifically, in the case of Haar multipliers, a pair of column/row operators are introduced by
\begin{equation*}
M_{\mathrm{c}}(f) = \sum_{Q \in \mathscr{D}} \alpha_Q \langle f , \phi_Q \rangle \phi_Q, \qquad M_{\mathrm{r}}(f) = \sum_{Q \in \mathscr{D}} \langle f , \phi_Q \rangle \alpha_Q  \phi_Q,
\end{equation*}
with uniformly bounded $\alpha_Q \in \mathcal{M}$. Even in the Lebesgue setting, Haar multipliers with noncommuting symbols may lack weak-type $(1,1)$ and strong $(p,p)$ estimates for $p \neq 2$.  This problem was solved in \cite{hong-lopezsanchez-martell-parcet2012} where weak-type $(1,1)$ estimates for Haar shift operators relative to the Lebesgue measure were obtained in terms of a column/row decomposition of the input function. To be more precise, given $f \in \mathcal{A}_{+,K}$ and $a, k \in \Z$ consider the Cuculescu's projections $q_k(2^c) = q_k(f,2^c)$ and
\begin{equation*}
\pi_{a,k} = \bigwedge_{c \geq a} q_k(2^c)\, - \bigwedge_{c \geq a - 1} q_k(2^c).
\end{equation*}
For fixed $k$ the projections $\pi_{a,k}$ are pairwise disjoint. Thus, $f$  decomposes in column/row components as $f = f_\mathrm{c} + f_\mathrm{r}$  in terms of the multiscale triangle truncations
\begin{equation*}
f_\mathrm{c} = \sum_{k\geq 1}\sum_{a \leq b} \pi_{a,k-1}df_k\pi_{b,k-1}\,,
\quad
f_\mathrm{r} = \sum_{k\geq 1}\sum_{a > b} \pi_{a,k-1}df_k\pi_{b,k-1}.
\end{equation*}
This decomposition is  used in \cite{hong-lopezsanchez-martell-parcet2012} in conjunction with the Calder\'{o}n-Zygmund decomposition found in \cite{parcet2009} to obtain that $\|M_\mathrm{r} f_\mathrm{r}\|_{1,\infty} + \|M_\mathrm{c} f_\mathrm{c} \|_{1,\infty} \, \lesssim \, \|f\|_1$, among analogous estimates for other Haar shift operators. Key to this argument is that the terms $\gamma$ in the Calder\'{o}n-Zygmund decomposition not having a proper $L_2$ estimate are such that
$\mathsf{D}_k (\gamma) = (1_{\mathcal{A}} - q_{k-1})A_k + A_k^*(1_{\mathcal{A}} - q_{k-1})$, which leads to vanishing  triangular truncations. A major setback for extending this argument to the nondoubling setting is that $\mathsf{D}_k(\beta_{\Delta}) = \beta_k = q_{k-1} \beta_k q_{k-1}$, reflecting that its classical counterpart decomposes into terms supported in the dyadic parents of the maximal dyadic cubes of $\Omega_\lambda$. This forces to estimate $L_1$ norms of triangular truncations of $\beta_k$, which in the $ \mathcal{B}(\ell_2^n)$-valued setting brings constants at best of order $\log (n+1)$. Furthermore, higher integrability of $\beta_k$ --- such as $L \log L$ (see \cite{randrianantoanina1998})--- might be hindered since $\mu$ is permitted to be nondoubling.

\end{remark}

\renewcommand\bibname{References}
\bibliographystyle{plain}
\small\bibliography{Refs}

\begin{thebibliography}{1}

\bibitem{condealoso-rey2014}
J.M. Conde-Alonso and G.~Rey.
\newblock On a pointwise estimate for positive dyadic shifts and some
  applications.
\newblock Preprint arXiv:1409.4351 [math.FA].

\bibitem{cuculescu1971}
I.~Cuculescu.
\newblock Martingales on von {N}eumann algebras.
\newblock {\em J. Multivariate Anal.}, 1(1):17--27, 1971.

\bibitem{hong-lopezsanchez-martell-parcet2012}
G.~Hong, L.D. L{\'o}pez-S{\'a}nchez, J.M. Martell, and J.~Parcet.
\newblock Calder\'on-{Z}ygmund operators associated to matrix-valued kernels.
\newblock {\em Int. Math. Res. Not.}, 2014(5):1221--1252, 2014.

\bibitem{hytonen2012}
T.~Hyt{\"o}nen.
\newblock The sharp weighted bound for general {C}alder\'on-{Z}ygmund
  operators.
\newblock {\em Ann. of Math. (2)}, 175(3):1473--1506, 2012.

\bibitem{lopezsanchez-martell-mparcet2014}
L.D. L\'{o}pez-S\'{a}nchez, J.M. Martell, and J.~Parcet.
\newblock Dyadic harmonic analysis beyond doubling measures.
\newblock {\em Adv. Math.}, 267:44--93, 2014.

\bibitem{mei-parcet2009}
T.~Mei and J.~Parcet.
\newblock Pseudo-lo\-ca\-li\-za\-tion of singular inte\-grals and
  non\-commuta\-tive {L}ittle\-wood\--{P}aley inequalities.
\newblock {\em Int. Math. Res. Not.}, 2009(8):1433--1487, 2009.

\bibitem{parcet2009}
J.~Parcet.
\newblock Pseudo-localization of singular integrals and noncommutative
  {C}alder\'on-{Z}ygmund theory.
\newblock {\em J. Funct. Anal.}, 256(2):509--593, 2009.

\bibitem{randrianantoanina1998}
N~Randrianantoanina.
\newblock Hilbert transform associated with finite maximal subdiagonal
  algebras.
\newblock {\em J. Austral. Math. Soc. Ser. A}, 65(3):388--404, 1998.

\bibitem{randrianantoanina2002}
N.~Randrianantoanina.
\newblock Non-commutative martingale transforms.
\newblock {\em J. Funct. Anal.}, 194(1):181--212, 2002.

\end{thebibliography}

\end{document}